\documentclass[12pt]{article}

\usepackage{amsmath,amssymb,amsthm}
 \def\d{\frac{1}{2}\,}
\def\sign {\mathrm{sign\, }}
\def\R{\mathbb{R}}
\def \<{\langle}\def\>{\rangle}

\newcommand{\ee}{\mathrm{e}} \newcommand{\dd}{\mathrm{d}} 
\renewcommand{\phi}{\varphi}

\title{The median of an exponential family and the normal law}

\author{G\'{e}rard Letac\thanks{Laboratoire de Statistique
et Probabilit\'es, Universit\'e  Paul Sabatier, Toulouse, France. \texttt{gerard.letac@math.univ-toulouse.fr}}, 
Lutz Mattner\thanks{Universit\"at Trier, Fachbereich IV - Mathematik, 54286 Trier, Deutschland.
\texttt{mattner@uni-trier.de}},
Mauro Piccioni\thanks{Dipartimento di Matematica, Sapienza Universit\`{a} di Roma, 00185 Roma, Italia.
\texttt{piccioni@mat.uniroma1.it}}}

=-0.5 in
\begin{document}
\maketitle 
\begin{abstract} Let $P$ be a probability on the real line  generating a natural exponential family  
$(P_t)_{t\in \R}$. We show that the property  that $t$ is  a median of $P_t$ for all $t$ characterizes 
$P$ as the standard Gaussian law $N(0,1).$


\vspace{2mm}
\end{abstract}

\section{Introduction}Let $P$ be a probability on the real line and  assume that

\begin{equation} \label{TLP}
L(t)=\int_{-\infty}^{+\infty}\ee^{tx}P(\dd x)<\infty\quad\text{ for }t\in\R.\end{equation}  Such a probability  generates the natural exponential family  $$\mathcal{F}_P=\{P_t(\dd x)=\frac{\ee^{tx}}{L(t)}P(\dd x),\  t\in \R  \}.$$ 
Then it might happen that the natural parameter $t$ of $\mathcal{F}_P$ is always a median of $P_t$, 
in the sense of 
\begin{equation}\label{PM}P_t((-\infty,t))\leq \d\leq P_t((-\infty,t]) 
\quad \text{ for }t\in\R.\end{equation}
In the sequel we denote by $\mathcal{P}$ the set of probabilities $P$ such that (\ref{TLP}) and (\ref{PM}) are fulfilled. A noteworthy example of an element of $\mathcal{P}$ is the standard normal distribution $N(0,1)$, for which $L(t)=\ee^{t^2/2}$ and $P_t=N(t,1).$  It will turn out that it is the only one. The following 
preliminary lemmas  simplify the study of $\mathcal{P}.$

\vspace{4mm}\noindent \textbf{Lemma 1}. {\em If $P\in \mathcal{P}$, then $P$ is absolutely continuous 
with respect to  Lebesgue measure. As a consequence, we have equality 
throughout in \eqref{PM}.} 

\vspace{4mm}\noindent \textbf{Lemma 2}.
{\em If $P\in \mathcal{P}$, then its distribution function is strictly increasing.}


\medskip
If $P\in \mathcal{P}$, then  Lemma 1 allows us to write
\begin{equation}\label{density}
P(\dd x)=g(x)\phi (x){\dd x},
\end{equation}
where $g$ is some 
measurable non-negative function and $\phi (x)=\ee^{-x^2/2}/\sqrt{2\pi}$ 
denotes the standard normal density, and we will show that then $g(x)=1$ a.e.~to get:

\vspace{4mm}\noindent \textbf{Theorem 1}. {\em 
$\mathcal{P}=\{N(0,1)\}$.}

\medskip
The proofs of the above results are contained in Section 2, followed by a conjecture and a further theorem.

\section{Proofs} 

\vspace{4mm}\noindent\textbf{Proof of Lemma 1}. 
The next paragraph shows that the distribution function of~$P$ is locally Lip\-schitz, and this implies 
the claimed absolute continuity, even with a locally bounded density, 
compare for example Royden and Fitzpatrick~(2010, pp.~120--124).

For $t\in\R$, multiplying in assumption \eqref{PM} by $L(t)$ yields
\begin{equation}    \label{Eq:h_versus_L}
  h(t):=\int_{(-\infty,t]}\ee^{tx}P(\dd x) \ge \frac{1}{2}L(t) \ge \int_{(-\infty,t)}\ee^{tx}P(\dd x) = h(t-).
\end{equation}
Hence, if $A>0$ is given, then for $s,t$ with $-A\le s < t \le A$, we get
\begin{eqnarray*}
 P\big( (s,t)\big) &=& \int_{(s,t)} \ee^{-tx}\ee^{tx}P(\dd x) \,\ \le  \,\ \ee^{A^2}\int_{(s,t)} \ee^{tx}P(\dd x) \\
  &=& \ee^{A^2} \left( h(t-) - h(s) + \int_{(-\infty,s]}(\ee^{sx}-\ee^{tx} )  P(\dd x) \right) \\
  &\le & \ee^{A^2} \left( \frac12(L(t)-L(s)) + (t-s)\int_\R |x| \ee^{A|x|}P(\dd x)\right) \\
  &\le & c^{}_A\cdot(t-s)
\end{eqnarray*}
for some finite constant $c^{}_A$.  We have been using ~\eqref{Eq:h_versus_L}
and $|\ee^u-\ee^v|\le |u-v|\ee^{w}$ for $|u|,|v|\le w$
at the penultimate step. Using  assumption~\eqref{TLP}, we rely     at the ultimate step
 on local Lipschitzness of $L,$  due to its analyticity, and on finiteness 
of $\int_\R |x| \ee^{A|x|}P(\dd x)$,  .  \qed

\vspace{4mm}\noindent \textbf{Proof of Lemma 2}. Assume to the contrary that 
there exist $a,b\in\R$ with $a<b$ and $P((a,b))=0$. 
Then, for  $t\in(a,b)$, Lemma~1 and \eqref{PM} yield 
 $$\int\limits_{-\infty}^a\ee^{tx}P(\dd x)
  =\int_{-\infty}^t \ee^{tx}P(\dd x)=\int_t^{+\infty} \ee^{tx}P(\dd x)=\int_b^{\infty}\ee^{tx}P(\dd x).$$ 
Thus the two measures $\pmb1_{(-\infty,a]}(x)P(\dd x)$ and $\pmb1_{[b,+\infty)}(x)P(\dd x)$ have 
finite and identical Laplace transforms on some non-empty interval. 
Hence the two measures coincide, and hence $P$ must be the zero measure, which is absurd. \qed

%

\vspace{4mm}\noindent \textbf{Proof of Theorem 1}. With the representation  (\ref{density}) for $P \in \mathcal {P}$, 
assumption \eqref{PM} is rewritten as 
\begin{equation} \label{Eq:BIE} 
\int_{-\infty}^t\ee^{tx-\frac{x^2}{2}}\frac{1}{\sqrt{2\pi}}g(x)\,\dd x=\frac 12 \int_{-\infty}^{+\infty}\ee^{tx-\frac{x^2}{2}}\frac{1}{\sqrt{2\pi}}g(x)\,\dd x.
\end{equation}
We multiply both sides by $\ee^{-t^2/2}:$ 
\begin{equation} \label{EIB}
\int_{-\infty}^t\ee^{-\frac{(t-x)^2}{2}}\frac{1}{\sqrt{2\pi}}g(x)\,\dd x=\frac 12 \int_{-\infty}^{+\infty}\ee^{-\frac{(t-x)^2}{2}}\frac{1}{\sqrt{2\pi}}g(x)\,\dd x.
\end{equation} 
In other terms the unknown function $g$ satisfies \begin{equation} \label{EQC}\int_{-\infty}^{+\infty}\sign(t-x)\phi(t-x)g(x)\,\dd x=0\end{equation} for all $t\in \R.$ A formal derivation of (\ref{EQC}) in $t$, 
using the product rule under the integral, and with one derivative being twice a delta function,
 leads to the equation
\begin{equation}\label{deriva} 
g(t)=\int_{-\infty}^{+\infty}  q(t-x) g(x) \,\dd x
\end{equation}
a.e.~in $t$, where $q(y):=\frac{1}{2}|y|\ee^{-\frac {y^2}{2}}$ is a probability density, 
but instead of justifying this formal differentiation, 
it seems easier to start by computing the derivative of 
$$
  h(t) := \int_{-\infty}^t \ee^{tx} P(\dd x).
$$

By Lemma~2 the distribution function $F$ of $P$ has a continuous inverse $F^{-1}$. 
Using the quantile transform we have
$$
  h(t) = \int_0^1 \pmb1^{}_{\{F^{-1}\le t\}}(u) \ee^{tF^{-1}(u)} \,\dd u  
    = \int_0^{F(t)} \ee^{tF^{-1}(u)} \,\dd u = H(F(t),t)
$$
with $H(s,t):= \int_0^s \ee^{tF^{-1}(u)} \,\dd u $ for $s\in(0,1)$ and $t\in\R$. Now $H$ has 
continuous partial derivatives  $H_1(s,t)= \ee^{tF^{-1}(s)} $ and $H_2(s,t)=\int_0^sF^{-1}(u) \ee^{tF^{-1}(u)} \,\dd u$,  
due to the continuity of  $F^{-1}$, 
and hence $H$ is differentiable. Let $f$ be a Lebesgue density of $P$. Then, at every $t$ where $F'(t)=f(t)$, and hence at Lebesgue-a.e.~$t$, the chain rule yields
\begin{eqnarray*}
 h'(t) &=& H_1(F(t),t)f(t) + H_2(F(t),t) \,\ = \,\ \ee^{t^2}f(t) + \int_0^{F(t)}F^{-1}(u) \ee^{tF^{-1}(u)} \,\dd u \\
  &=&  \ee^{t^2}f(t) + \int_{-\infty}^t x\ee^{tx}f(x)\,\dd x.
\end{eqnarray*}
Thus differentiating the identity~\eqref{Eq:BIE} and observing that $f(x)=g(x)\phi (x)$ we obtain the 
following a.e.-identity 
$$
  \frac{1}{\sqrt{2\pi}}\ee^{t^2/2}g(t) + \int_{-\infty}^t x \ee^{tx-\frac{x^2}2}\frac{1}{\sqrt{2\pi}}g(x)\,\dd x 
= \frac12\int_{-\infty}^{+\infty} x \ee^{tx-\frac{x^2}2}\frac{1}{\sqrt{2\pi}}g(x)\,\dd x,
$$
and multiplying the latter by $\sqrt{2\pi}\ee^{-t^2/2}$ gives
$$
g(t)=\frac 12 \biggl(\int_t^{+\infty} x\ee^{-(t-x)^2/2}g(x)\,\dd x-\int_{-\infty}^t x\ee^{-(t-x)^2/2}g(x)\,\dd x \biggr).
$$
Adding to the rigth hand side above the quantity 
$$
0=\frac {t}{2}\biggl(\int_{-\infty}^t \ee^{-(t-x)^2/2}g(x)\,\dd x-\int_t^{+\infty}\ee^{-(t-x)^2/2}g(x)\,\dd x\biggr)
$$
(recall (\ref{EIB})) yields the desired~\eqref{deriva}.

Next, with the (positive) Radon measures $\mu(\dd x):=g(x)\dd x$ and $\sigma(\dd x):=q(x)\dd x$,
equation~\eqref{deriva} can be rewritten as the so-called Choquet-Deny equation  $\mu= \mu \ast \sigma$.  
Observe that  $t\mapsto \int_{-\infty}^{+\infty} \ee^{tx}\sigma (\dd x)$ is even and strictly convex,
and is therefore  equal to~$1$ only at $t=0.$ 
We can now use the results in section 6 of Deny~(1960), where ``$n>1$'' is evidently a misprint 
for ``$n\ge1$'', to conclude that 
$\mu$ has to be a positive scalar multiple of the Lebesgue measure. Since $g$ is a probability 
density with respect to a probability measure, we have $g=1$ a.e., and the theorem is proved. \qed

\medskip
Finally, it is worthwhile to mention a natural conjecture about exponential families which seems harder to establish: 

\vspace{4mm}\noindent\textbf{Conjecture.} Suppose  that the probability $P$ satisfies (\ref{TLP}), and denote  $m(t):=\int_{\R}xP_t(\dd x).$  If for all $t$ real $m(t)$ is a median of $P_t$, then $P=N(m,\sigma^2)$ for some $m$ and $\sigma.$

\medskip
This conjecture, which is probably more meaningful from a methodological point of view than the result established in the paper, does not translate in a neat harmonic analysis statement as (\ref{EQC}) and (\ref{deriva}) and as such it seems harder to establish. The next simple result offers some support to the conjecture. A probability $Q$ on $\R^n$ is said to be symmetric  if there exists some  $m\in \R^n$ such that $X-m\sim m-X$ when $X\sim Q.$

\vspace{4mm}\noindent\textbf{Theorem 2.} {\em Let $P$ be a probability on $\R^n$ such that 
$$L(t)=\int _{\R^n}\ee^{\<t,x\>}P(\dd x)$$ is finite for all $t\in \R^n.$ Assume that for all $t\in \R^n$ the probability $P_t(\dd x)=\ee^{\<t,x\>}P(\dd x)/L(t)$ is symmetric. Then $P$ is normal.}

\vspace{4mm}\noindent\textbf{Proof.} Clearly $m(t)=\int_{\R^n}xP_t(\dd x)=L'(t)/L(t)$ exists and, since $P_t$ is symmetric, $X_t-m(t)\sim m(t)-X_t$ when $X_t\sim P_t.$ Therefore its Laplace transform 
$$s\mapsto \mathbb{E}(\ee^{\<s,X_t-m(t)\>})=\ee^{-\<s,m(t)}\frac{L(t+s)}{L(t)}$$
does not change when we replace $s$ by $-s.$ Considering  the logarithm and taking the derivative in $s$ we get $2m(t)=m(t+s)+m(t-s)$. Taking again the derivative in $s$ we get $m'(t+s)=m'(t-s)$ for all $t,s\in \R^n$, which means that $m'$ is constant, hence $\log L$ is polynomial of degree at most 2,
 and hence $P$ is normal. \qed

\section{References}
\vspace{4mm}\noindent \textsc{Deny, J.} (1960). Sur l'\'equation de convolution $\mu=\mu \ast \sigma$. 
\textit{S\'eminaire Brelot-Choquet-Deny (Th\'eorie du Potentiel)} \textbf{4}e ann\'ee, 1959-60,  
Expos\'e num\'ero 5.\\ \verb{http://www.numdam.org/article/SBCD_1959-1960__4__A5_0.pdf{

\vspace{4mm}\noindent \textsc{Royden, H.L.} and \textsc{Fitzpatrick, P.M} (2010).
\textit{Real Analysis.} Fourth edition, Prentice-Hall.

\end{document}